\documentclass[12pt,a4paper]{amsart}
\usepackage[utf8]{inputenc}	
\usepackage[T1]{fontenc}
\usepackage[in,headings]{fullpage}
\usepackage{amsbsy, amscd, amsfonts, amsmath, amsrefs, amssymb, amsthm}
\usepackage{graphicx}
\usepackage{float}
\usepackage{color}   
\usepackage[colorlinks=true,linkcolor=blue,citecolor=blue,urlcolor=blue]{hyperref}

\newtheorem{theorem}{Theorem}%[section]

\newtheorem{corollary}[theorem]{Corollary}
\newtheorem{proposition}[theorem]{Proposition}
\newtheorem{lemma}[theorem]{Lemma}

\newtheorem*{conjecture}{Conjecture}

\theoremstyle{definition}

\newtheorem{remark}[theorem]{Remark}

\theoremstyle{remark}

\newcommand{\R}{\mathbf{R}}

\renewcommand{\Re}{\mathop{\mathrm{Re}}\nolimits}
\renewcommand{\Im}{\mathop{\mathrm{Im}}\nolimits}
\newcommand{\Rzeta}{\mathop{\mathcal R }\nolimits}

\newfont{\cmbsy}{cmbsy10}
\newfont{\cmmib}{cmmib10}

\begin{document}

\title{Riemann's Auxiliary Function. Right limit of zeros}
\author[Arias de Reyna]{J. Arias de Reyna}
\address{%
Universidad de Sevilla \\ 
Facultad de Matem\'aticas \\ 
c/Tarfia, sn \\ 
41012-Sevilla \\ 
Spain.} 

% AMS subject classifications (used in AMS journals)

% AMS keywords (used in AMS journals)

% acknowledge support, etc
%\thanks{This research was  supported by MINECO grant MTM2015--63699-P}
% \thanks{We would like to thank our colleagues for their helpful
%  criticism.}

\email{arias@us.es, ariasdereyna1947@gmail.com}

%\date{\today, \texttt{173-RightLimit-v6.tex}}

\begin{abstract}
Numerical data suggest that the zeros $\rho$ of the auxiliary Riemann function in the upper half-plane satisfy $\Re(\rho)<1$. We show that this is true for those zeros with $\Im(\rho)> 3.9211\dots10^{65}$. We conjecture that this is true for all of them.
\end{abstract}

\maketitle

\section{Introduction}

We denote by $\Rzeta(s)$ the auxiliary  Riemann function (see \cite{A166} for its definition), which, according to Siegel \cite{Siegel}, was considered by Riemann in his Nachlass.
In this paper, we consider the zeros of $\Rzeta(s)$ in the first quadrant. Siegel related the zeros of $\Rzeta(s)$  to the zeros of $\zeta(s)$. Each zero $\rho=\beta+i\gamma$ of $\Rzeta(s)$,  with $\beta<1/2$ and $\gamma>0$ generates two  $\zeta(s)$ zeros on the critical line. But here we are considering subsets of the first quadrant without zeros of $\Rzeta(s)$. 

In \cite{A172} I give some statistics based on my computation of the  zeros in the domain
$0<t<200000$. All these zeros $\rho=\beta+i\gamma$ satisfies $\beta<1$. My data suggest that 
$\sup_{\gamma>0}\beta=1$.  Therefore, we may state the conjecture:
\begin{conjecture}
Any zero $\rho=\beta+i\gamma$ of $\Rzeta(s)$, with $\Im(\rho)=\gamma>0$ satisfies $\beta<1$. 
\end{conjecture}
This paper contains my results on this conjecture. I prove it for zeros with $\gamma\ge t_0$, but for a $t_0$ too large.

Siegel \cite{Siegel}*{eq.(84)} was the first to obtain some positive result in this direction. He proved that $|\Rzeta(s)-1|<3/4$ for $\sigma\ge2$ and $t>t_1$ for some unspecified $t_1$. Here we apply the asymptotic development \cite{A86} of  $\Rzeta(s)$ to make Siegel's results more precise.  
\begin{proposition} 
For $\sigma\ge 2$ and $t\ge 32\pi$, we have $|\Rzeta(s)-1|<1$, so that $\Rzeta(s)\ne0$ in 
$[2,+\infty)\times[32\pi,+\infty)$.
\end{proposition}

In \cite{A100}*{Corollary 14} it is proved that there exists $r_0>0$ such that for $\sigma\ge2$, $t\ge0$ and $|s|>r_0$ we have $|\Rzeta(s)-1|<3/4$. Also, my computation of zeros of $\Rzeta(s)$ make very improbable any zero with $\gamma>0$ and $\beta\ge2$.

Using the same expansion of $\Rzeta(s)$ and  Rouche's Theorem (see Ahlfors \cite{Ah}*{p.~153}) we can extend the region without zeros:
\begin{proposition}
The function $\Rzeta(s)$ does not vanish in the rectangle $[3/2,2]\times[2707,+\infty)$.
\end{proposition}

Finally, adding van der Corput estimates, we get our main result.
\begin{theorem}\label{T:main}
The auxiliary function $\Rzeta(s)$ do not vanish for $1\le \sigma\le2$ and $t>t_0=3.9211\dots10^{65}$.
\end{theorem}
We would like to decrease  the value of $t_0$. To this end we   need either a better bound of $|\zeta(\sigma+it)|^{-1}$ for $\sigma\approx 1$ than that given by  Trudgian \cite{Tr2} and Mossinghoff and Trudgian  \cite{MT}, or a better bound of $\zeta(s)-\sum_{n<\sqrt{t/2\pi}}n^{-\sigma-it}$ for $\sigma\ge1$  than that given in Proposition \ref{P:vanderCorput}.

We prove Theorem \ref{T:main} applying Rouche's theorem \cite{Ah}*{p.~153}. To this end, we show (see Proposition \ref{P:RzetaApp}) some interesting inequalities for the auxiliary function on the critical strip and bound the difference $\zeta(s)-\sum_{n\le\sqrt{t/2\pi}}n^{-\sigma-it}$ for $\sigma>1$.

\section{Approximation of \texorpdfstring{$\Rzeta(s)$}{R(s)} by zeta sums.}

\begin{proposition}\label{P:RzetaApp}
The following inequalities hold:
\begin{gather}\label{E:boundRzeta}
\Bigl|\Rzeta(s)-\sum_{n\le \sqrt{t/2 \pi}}\frac{1}{n^s}\Bigr|\le \Bigl(\frac{t}{2\pi}\Bigr)^{-\sigma/2},\qquad \text{for $0\le\sigma\le1$ and $t\ge3\pi$},\\
\Bigl|\Rzeta(s)-\sum_{n\le \sqrt{t/2 \pi}}\frac{1}{n^s}\Bigr|\le \Bigl(\frac{t}{2\pi}\Bigr)^{-\sigma/2},\qquad \text{for $1\le\sigma\le2$ and $t\ge 8\pi$},\\
\Bigl|\Rzeta(s)-\sum_{n\le \sqrt{t/2 \pi}}\frac{1}{n^s}\Bigr|\le \Bigl(\frac{t}{2\pi}\Bigr)^{-1/2},\qquad \text{for $1\le\sigma$ and $t\ge 16\pi$}.
\end{gather}
\end{proposition}

\begin{proof}
Given $s=\sigma+it$ with  $\sigma>0$ and $t>0$ and taking $K=1$ in
\cite{A86}*{Theorem 3.1, 4.1 and 4.2} we obtain with $a=\sqrt{t/2\pi}$, $N=\lfloor a\rfloor$, $p=1-2(a-N)$ and $|U|=1$ 
\[\Rzeta(s)=\sum_{n\le a}\frac{1}{n^s}+(-1)^{N-1} U a^{-\sigma}\Bigl\{C_0(p)+\frac{C_1(p)}{a}+RS_1\Bigr\}.\]
We also have the following bounds
\[\frac{|C_1(p)|}{a}\le \frac{\sqrt{2}}{2\pi} \frac{9^\sigma\sqrt{\pi}}{2 a}, \qquad
|RS_1|\le \frac{2^{3\sigma/2}}{7}\frac{1}{(10a/11)^2}.\]
Also in \cite{A86}*{eq. (5.2) and (6.7)} we find $|C_0(p)|\le\frac12$. 

These bounds, in particular those for $C_j(p)$, are large,and can be improved if we use the values given in \cite{A86}*{(2.4)}
\begin{align*}
C_0(p)&=F(p),\\
C_1(p)&=\frac{1}{\pi^2}\Bigl(d^{(1)}_0 F^{(3)}(p)+\frac{\pi}{2i}d^{(1)}_1F'(p)\Bigr)=\frac{i(\sigma-\frac12)}{2\pi}F'(p)+\frac{1}{12\pi^2}F^{(3)}(p),
\end{align*}
Since the recurrence equation for the coefficients $d^{(k)}_j$
\cite{A86}*{(2.7), (2.8), (2.9) and (2.10)} gives us $d^{(1)}_0=\frac{1}{12} $ and $d^{(1)}_1=\frac12-\sigma$

Combining these with the bounds for $F^{(n)}(p)$ given in \cite{A86}*{(6.7)}
\[|F^{(2n)}(p)|\le\frac{(2n)!}{2^{n+1}n!}\pi^n, \qquad F^{(2n+1)}(p) \le 2^n\pi^n n!\]
yields
\[|C_1(p)|\le \frac{|\sigma-\frac12|}{2\pi}+\frac{1}{6\pi}\]
Therefore,
\begin{equation}\label{E:partial}
\Bigl|\Rzeta(s)-\sum_{n\le \sqrt{t/2 \pi}}\frac{1}{n^s}\Bigr|\le a^{-\sigma}\Bigl(\frac12+\frac{1}{6\pi a}+\frac{|\sigma-\frac12|}{2\pi a}+\frac{2^{3\sigma/2}}{7}\frac{1}{(10a/11)^2}\Bigr).
\end{equation}
For $0\le\sigma\le1$ we get 
\[\Bigl|\Rzeta(s)-\sum_{n\le \sqrt{t/2 \pi}}\frac{1}{n^s}\Bigr|\le a^{-\sigma}\Bigl(\frac12+\frac{1}{6\pi a}+\frac{1}{4\pi a}+\frac{2^{3\sigma/2}}{7}\frac{1}{(10a/11)^2}\Bigr).\]
This is less than $a^{-\sigma}$ for $a\ge\sqrt{3/2}$, and hence for
$t\ge3\pi$.

For $\sigma\ge1$ and $a>1$ we have
\[\Bigl|\Rzeta(s)-\sum_{n\le \sqrt{t/2 \pi}}\frac{1}{n^s}\Bigr|\le \frac{1}{2a}+\frac{1}{6\pi a^2}+\frac{|\sigma-\frac12|}{2\pi a\cdot a^\sigma}+\frac{1}{7}(2^{3/2}/a)^\sigma\frac{1}{(10a/11)^2}.\]
When $t\ge16\pi$, we get $2^{3/2}/a \le 1$ and $0<2\sigma-1<2\sigma<2^{3\sigma/2}\le a^\sigma$
and 
\[\Bigl|\Rzeta(s)-\sum_{n\le \sqrt{t/2 \pi}}\frac{1}{n^s}\Bigr|\le \frac{1}{a}\Bigl(\frac{1}{2}+\frac{1}{6\pi a}+\frac{1}{4\pi}+\frac{1.21}{7a}\Bigr),\]
and this is less than $1/a$ when $a>1$.

For $1\le\sigma\le 2$, \eqref{E:partial} gives us
\[\Bigl|\Rzeta(s)-\sum_{n\le \sqrt{t/2 \pi}}\frac{1}{n^s}\Bigr|\le a^{-\sigma}\Bigl(\frac12+\frac{1}{6\pi a}+\frac{3}{4\pi a}+\frac{8}{7}\frac{1}{(10a/11)^2}\Bigr).\]
This is less than $a^{-\sigma}$ for $a\ge2$, hence for $t\ge8\pi$. 
\end{proof}

\begin{remark}
The inequalities in Proposition \ref{P:RzetaApp} can be improved by taking more terms of the Riemann-Siegel expansion, bounding the terms in a way similar to what we have done here.
\end{remark}

\begin{proposition} For $\sigma\ge 2$ and $t>16\pi$ we have
\begin{equation}
|\Rzeta(s)-1|<
\frac{3}{2^\sigma}+\Bigl(\frac{t}{2\pi}\Bigr)^{-1/2}.
\end{equation}
Hence, for $\sigma\ge 2$ and $t\ge32\pi$ we have $|\Rzeta(s)-1|<1$, and there is no zero in this closed set.
\end{proposition}

\begin{proof} By Proposition \ref{P:RzetaApp} for $\sigma\ge2$ and $t\ge 16\pi$ we have
\begin{align*}
|\Rzeta(s)-1|&\le \sum_{n=2}^\infty
\frac{1}{n^\sigma}+\Bigl(\frac{t}{2\pi}\Bigr)^{-1/2}\le \frac{1}{2^\sigma}+\frac{1}{3^\sigma}+
\int_3^{+\infty}\frac{dv}{v^\sigma}+\Bigl(\frac{t}{2\pi}\Bigr)^{-1/2}\\
&\le \frac{1}{2^\sigma}+\frac{4}{3^\sigma}
+\Bigl(\frac{t}{2\pi}\Bigr)^{-1/2},
\end{align*}
and since $\sigma\ge 2$
\[<\frac{3}{2^\sigma}+\Bigl(\frac{t}{2\pi}\Bigr)^{-1/2}.\]
When $\sigma\ge2$ and $t\ge 32\pi$ we have 
\[|\Rzeta(s)-1|<\frac{3}{2^\sigma}+\Bigl(\frac{t}{2\pi}\Bigr)^{-1/2}\le \frac{3}{4}+\frac14=1.\qedhere\]
\end{proof}

\section{First application of Rouche's Theorem}\label{S:3}

In this section, we show that there is no zero of $\Rzeta(s)$ in the rectangle $[3/2,2]\times[t_0.+\infty)$. In this range of $\sigma$ the result is almost trivial. We have to prove 
$|\Rzeta(s)-\zeta(s)|<|\zeta(s)|$ for $s$ on the boundary of the region. We need a lower bound for $\zeta(s)$ given by the next Lemma.
\begin{lemma}\label{L:1}
For $\sigma>1$ we have 
\begin{equation}
|\zeta(s)|\ge \frac{\zeta(2\sigma)}{\zeta(\sigma)}.
\end{equation}
\end{lemma}
\begin{proof}
The proof is simple 
\[\Bigl|\frac{1}{\zeta(s)}\Bigr|=\Bigl|\prod_p\Bigl(1-\frac{1}{p^s}\Bigr)\Bigr|\le 
\prod_p\Bigl(1+\frac{1}{p^\sigma}\Bigr)=\frac{\zeta(\sigma)}{\zeta(2\sigma)}.\qedhere\]
\end{proof}

\begin{proposition}
The function $\Rzeta(s)$ does not vanish in the rectangle $[3/2,2]\times[2707,+\infty)$
\end{proposition}
\begin{proof}
Take $T>t_0=2707$. We apply Rouche's Theorem to the rectangle $R=[3/2,2]\times[t_0,T]$. We need to  prove the inequality
\[|\Rzeta(s)-\zeta(s)|< |\zeta(s)|,\qquad\text{for $s$ in the boundary of $R$}.\]
We may start by using Proposition \ref{P:RzetaApp}. 
\[|\Rzeta(s)-\zeta(s)|\le \Bigl|\Rzeta(s)-\sum_{n\le \sqrt{t/2 \pi}}\frac{1}{n^s}\Bigr|+
\Bigl|\sum_{n\le \sqrt{t/2 \pi}}\frac{1}{n^s}-\zeta(s)\Bigr|.\]
Since $t_0>8\pi$,  for $\sigma\in[3/2,2]$ we have
\[|\Rzeta(s)-\zeta(s)|\le\Bigl(\frac{t}{2\pi}\Bigr)^{-\frac{\sigma}{2}}+\Bigl|\sum_{n>\sqrt{t/2 \pi}}\frac{1}{n^\sigma}\Bigr|.\]
Substituting the sum  by a first term and an integral, we get
\[|\Rzeta(s)-\zeta(s)|\le2\Bigl(\frac{t}{2\pi}\Bigr)^{-\frac{\sigma}{2}}+\frac{1}{\sigma-1}
\Bigl(\frac{t}{2\pi}\Bigr)^{\frac{1-\sigma}{2}}.\]

For each $\sigma>1$ there is a $t_0$ such that for $t>t_0$
\[2\Bigl(\frac{t}{2\pi}\Bigr)^{-\frac{\sigma}{2}}+\frac{1}{\sigma-1}
\Bigl(\frac{t}{2\pi}\Bigr)^{\frac{1-\sigma}{2}}< \frac{\zeta(2\sigma)}{\zeta(\sigma)}\le|\zeta(s)|.\]
This will give us the inequality $|\Rzeta(s)-\zeta(s)|< \zeta(s)$.  The function $\frac{\zeta(2\sigma)}{\zeta(\sigma)}$ is increasing with $\sigma$ as we see from the proof of Lemma \ref{L:1}, therefore, for $\frac32\le\sigma\le 2$ we have for $t\ge2707$ 
\begin{multline*}
2\Bigl(\frac{t}{2\pi}\Bigr)^{-\frac{\sigma}{2}}+\frac{1}{\sigma-1}
\Bigl(\frac{t}{2\pi}\Bigr)^{\frac{1-\sigma}{2}}\le 2\Bigl(\frac{t}{2\pi}\Bigr)^{-\frac{\sigma}{2}}+2
\Bigl(\frac{t}{2\pi}\Bigr)^{\frac{1-\sigma}{2}}\\
\le 2\Bigl(\frac{t}{2\pi}\Bigr)^{-\frac34}+2
\Bigl(\frac{t}{2\pi}\Bigr)^{-\frac{1}{4}}
< \frac{\zeta(3)}{\zeta(3/2)} \le\frac{\zeta(2\sigma)}{\zeta(\sigma)}\le|\zeta(s)|,
\end{multline*}
where the strict inequality is checked by numerically computing the left and right side for $t=2707$.
Applying Rouche's theorem, we will find that $\Rzeta(s)$ and $\zeta(s)$ have the same number of zeros in the rectangle $R$. Since $\zeta(s)$ do not vanish, we get that $\Rzeta(s)\ne0$ for $s\in R$. Since $T>t_0$ is arbitrary, the proposition is proved.
\end{proof}

We want now to apply Rouché's Theorem for a rectangle of type $[1,3/2]\times[t_0,T]$.  For $\sigma=1$ the Dirichlet series for $\zeta(s)$ do not apply, and the argument will be more complicated. 

In particular, we will need  explicit bounds of 
\[\Bigl|\zeta(s)-\sum_{n\le \sqrt{t/2 \pi}}\frac{1}{n^s}\Bigr|,\qquad 1\le\sigma\le 3/2,\quad t\ge t_0,\]
and explicit bounds in the same range of $1/\zeta(s)$.

\section{Approximating \texorpdfstring{$\zeta(s)$}{zeta(s)} by a zeta sum}

\subsection{Approximating \texorpdfstring{$\zeta(s)$}{zeta(s)} by a long zeta sum.}
For $\sigma>0$ and any real number $x>0$ by partial summation, we get 
\begin{equation}\label{E:EulerMac}
\zeta(s)=\sum_{n\le x}\frac{1}{n^s}+\frac{x^{1-s}}{s-1}+\frac{\{x\}-\frac12}{x^s}+s\int_x^\infty
\frac{\frac12-\{u\}}{u^{s+1}}\,du.
\end{equation}
As a matter of notation, we will use $\tau=\frac{t}{2\pi}$ as an abbreviation. We always consider $s=\sigma+it$ with $\sigma\in\R$ and $t>0$. 
\begin{proposition}\label{P:MacLaurin}
For $\sigma\ge1/2$, $\tau=\frac{t}{2\pi}>1$ and $0<r\le 2$
\begin{equation}
\Bigl|\zeta(s)-\sum_{n\le \tau^r}\frac{1}{n^s}\Bigr|\le  8\tau^{1-r\sigma}.
\end{equation}
\end{proposition}
\begin{proof}
In \eqref{E:EulerMac} put $x=\tau^r$, and notice that
\[\Bigl|\frac{\tau^{r(1-s)}}{s-1}\Bigr|\le \frac{\tau^{r(1-\sigma)}}{|\sigma-1+it|}\le \frac{1}{2\pi}
\tau^{r(1-\sigma)-1},\quad \Bigl|\frac{\{\tau^r\}-\frac12}{(\tau^r)^s}\Bigr|\le \frac{1}{2\tau^{r\sigma}},\]
\[\Bigl|s\int_{\tau^r}^\infty \frac{1/2-\{u\}}{u^{s+1}}\,du\Bigr|\le 
\frac{|\sigma+it|}{2}\int_{\tau^r}^\infty\frac{du}{u^{\sigma+1}}=\frac{|\sigma+it|}{2}\frac{\tau^{-r\sigma}}{\sigma}= \pi\Bigl|\frac{1}{t}+\frac{i}{\sigma}\Bigr|\tau^{1-r\sigma}.
\]
Hence,
\begin{align*}
\Bigl|\zeta(s)-\sum_{n\le \tau^r}\frac{1}{n^s}\Bigr|&\le \frac{1}{2\pi}
\tau^{r(1-\sigma)-1}+\frac12\tau^{-r\sigma}+\pi\Bigl(\frac{1}{t}+\frac{1}{\sigma}\Bigr)\tau^{1-r\sigma}\\
&\le \Bigl(\frac{1}{2\pi\tau^{2-r}}+\frac{1}{2\tau}+\frac{\pi}{t}+2\pi\Bigr)\tau^{1-r\sigma}\\
&\le \Bigl(\frac{1}{2\pi}+1+2\pi\Bigr)\tau^{1-r\sigma}\le 8\tau^{1-r\sigma}.\qedhere
\end{align*}
\end{proof}

\begin{remark}
Kadiri \cite{K} gives an improvement of Proposition \ref{P:MacLaurin}. But this improvement gives us only  a slightly better $t_0$ in our main Theorem \ref{T:main}.
\end{remark}

\subsection{van der Corput's machinery.}
We need to approximate $\zeta(s)$ by the sum $\sum_{n\le \tau^{1/2}}\frac{1}{n^{s}}$ for $1\le\sigma\le 3/2$. We will use the van der Corput $d$-th derivative test, for which we have an explicit form.
Although these explicit forms can be found, for example, in Platt and Trudgian \cite{PT} and   Hiary \cite{H}, we do not use these versions since they are founded on a wrong lemma.  Instead, we use \cite{A181}.

For $X>0$, $t>1$ and $\sigma>0$ we define as in \cite{A181}  the bounds of the zeta sums
\begin{equation}
S(X,t):=\sup_{X<Z\le 2X}\Bigl|\sum_{X<n\le Z}n^{-it}\Bigr|,\qquad S_\sigma(X,t):=\sup_{X<Z\le 2X}\Bigl|\sum_{X<n\le Z}\frac{1}{n^{\sigma+it}}\Bigr|.
\end{equation}
As a matter of notation, it is convenient to put $\tau=t/2\pi$. 
Define the exponent $\alpha$ of $S(X,t)$ as the real number such that $X=\tau^\alpha$. We will use $\alpha$ with this meaning in the bounds of $S(X,t)$. 
By partial summation, see \cite{A181}*{Prop.~22} we see 
\begin{equation}\label{E:Abel}
S_\sigma(X,t)\le X^{-\sigma}S(X,t).
\end{equation}

We will need  Propositions 18 and 20  of \cite{A181}:

\begin{proposition}[Th. 18 in \cite{A181}]\label{P:d=2}
Let $\alpha>0$ and $t>2\pi$, then 
\begin{equation}\label{E:191218-1}
S(\tau^\alpha,t)\le 2A(\tau^{1/2}+2\tau^{\alpha-1/2}).
\end{equation}
\end{proposition}

and 

\begin{proposition}[Th. 20 in \cite{A181}]\label{L:Sbound}
Let  $\alpha>0$ be a real number and $d\ge2$  an integer.  Assume that  $\alpha< \frac{2^{d-2}}{1+(d-2)2^{d-2}}$. For any $M\ge \widehat{B_d}$,  there exists a number $\tau_0=\tau_0(M)$ such that
\[S(\tau^\alpha,t)\le M \tau^{\alpha+\frac{1-\alpha d}{D-2}}\qquad \tau=t/2\pi\ge\tau_0.\]
The number $\tau_0$ is the least number such that  
\begin{equation}\label{E:t0}
\tau_0\ge d^{D/\alpha(D-2)},\quad \tau_0^{\frac{2}{D}+\frac{1-\alpha d}{D-2}}\ge\widehat{C_d}/M,\quad \tau_0^{\frac{2\alpha}{D}+\frac{1-\alpha d}{D-2}}\ge\widehat{A_d}/M.\end{equation}
\end{proposition}
Here $A$, $B$, $\widehat{A_d}$, $\widehat{B_d}$, and $\widehat{C_d}$ are constants defined in \cite{A181} and $D=2^d$. 

\subsection{Approximating \texorpdfstring{$\zeta(s)$}{zeta(s)} by a short zeta sum.}
\begin{proposition}\label{P:van1}
For $\frac12< \sigma< 1$ and $\tau\ge 20$, we have
\[\Bigl|\zeta(s)-\sum_{n\le \tau^{1/2}}\frac{1}{n^s}\Bigr|\le 
\frac{\widehat{A_3}2^{\sigma-\frac12}}{2^{\sigma-\frac12}-1}\tau^{\frac12(\frac56-\sigma)}+
\frac{2^{1+\sigma}A}{1-2^{-\sigma}}\tau^{\frac23(\frac34-\sigma)}+4A\Bigl(1+\frac{4\log\tau}{3\log2}\Bigr)\tau^{2(\frac34-\sigma)}+8\tau^{2(\frac12-\sigma)}.\]
\end{proposition}
\begin{proof}
By Proposition \ref{P:MacLaurin} with $r=2$ 
\[\Bigl|\zeta(s)-\sum_{n\le \tau^{1/2}}\frac{1}{n^s}\Bigr|\le 8\tau^{1-2\sigma}+\Bigl|\sum_{\tau^{1/2}<n\le \tau^{2}}\frac{1}{n^s}\Bigr|.\]
To bound the last sum, notice that the exponents of our sums will be $\frac12\le \alpha\le 2$. As indicated in Remark 19 in \cite{A181}, we must use the van der Corput $d$-th derivative test for $d=2$ or $d=3$ according to whether it is $\alpha\ge 2/3$ or not. 
Hence, we determine  a natural number $K$ such that $2^{-K}\tau^2<\tau^{2/3}\le 2^{-K+1}\tau^2$ and split the sum 
\[\Bigl|\sum_{\tau^{1/2}<n\le\tau^2}\frac{1}{n^{\sigma+it}}\Bigr|\le
\Bigl|\sum_{\tau^{1/2}<n\le2^{-K}\tau^2}\frac{1}{n^{\sigma+it}}\Bigr|+
\Bigl|\sum_{2^{-K}\tau^{2}<n\le\tau^2}\frac{1}{n^{\sigma+it}}\Bigr|.\]
Now there is some integer $L\ge0$ with $2^L\tau^{1/2}< 2^{-K}\tau^{2}\le 2^{L+1}\tau^{1/2}$. We have the following
\[\Bigl|\sum_{\tau^{1/2}<n\le\tau^2}\frac{1}{n^{\sigma+it}}\Bigr|\le
\sum_{\ell=0}^L\Bigl|\sum_{2^\ell\tau^{1/2}<n\le2^{\ell+1}\tau^{1/2}}\frac{1}{n^{\sigma+it}}\Bigr|
+\sum_{k=1}^K \Bigl|\sum_{2^{-k}\tau^{2}<n\le2^{-k+1}\tau^2}\frac{1}{n^{\sigma+it}}\Bigr|,\]
where in the sum corresponding to $\ell=L$ with the range $2^L\tau^{1/2}<n\le2^{L+1}\tau^{1/2}$ we only include
the terms corresponding to the numbers $n$ satisfying $2^L\tau^{1/2}<n\le 2^{-K}\tau^2\le 2^{L+1}\tau^{1/2}$.
Therefore, by \eqref{E:Abel} and the definition of the $S(X,t)$
\[\Bigl|\sum_{\tau^{1/2}<n\le\tau^2}\frac{1}{n^{\sigma+it}}\Bigr|\le
\sum_{\ell=0}^L(2^\ell\tau^{1/2})^{-\sigma}S(2^\ell\tau^{1/2},t)
+\sum_{k=1}^K (2^{-k}\tau^{2})^{-\sigma}S(2^{-k}\tau^{2},t):= S_1+S_2.\]
To the second sum, we apply Proposition \ref{P:d=2}, obtaining 
\[S_2=\sum_{k=1}^K (2^{-k}\tau^{2})^{-\sigma}S(2^{-k}\tau^{2},t)\le 2A\sum_{k=1}^K
(2^{-k}\tau^{2})^{-\sigma}\bigl\{\tau^{1/2}+2\tau^{\alpha_k-1/2}\bigr\},\]
where $\alpha_k$  is the exponent such that $\tau^{\alpha_k}=2^{-k}\tau^{2}$. Therefore,
\begin{equation}\label{E:partial1}
S_2\le \sum_{k=1}^K (2^{-k}\tau^{2})^{-\sigma}S(2^{-k}\tau^{2},t)\le 2A\sum_{k=1}^K
(2^{-k}\tau^{2})^{-\sigma}\bigl\{\tau^{1/2}+2\tau^{-1/2}2^{-k}\tau^{2}\}.
\end{equation}
Since we assume that $\frac12\le \sigma<1$, we continue with
\[S_2\le2A\tau^{\frac12-2\sigma}\frac{2^{K\sigma}-1}{1-2^{-\sigma}}+4A\tau^{\frac32-2\sigma}\frac{1-2^{-K(1-\sigma)}}{2^{1-\sigma}-1}.\]
When $\sigma$ is near $1$ the second expression is not very useful. We have 
$\frac{1-2^{-K(1-\sigma)}}{2^{1-\sigma}-1}=\sum_{k=1}^K 2^{-k(1-\sigma)}\le K$. Hence, 
this second sum is bounded by 
\[S_2\le \frac{2A}{1-2^{-\sigma}}2^{K\sigma}\tau^{\frac12-2\sigma}+4AK\tau^{\frac32-2\sigma}.\]
By its definition $2^K\le 2\tau^{4/3}$ so that 
\[S_2\le \frac{2^{1+\sigma}A}{1-2^{-\sigma}}\tau^{\frac12-\frac{2}{3}\sigma}+4A\Bigl(1+\frac{4\log\tau}{3\log2}\Bigr)\tau^{\frac32-2\sigma}.\]
To the first sum, we apply Proposition \ref{L:Sbound} with $d=3$. For $S(2^\ell\tau^{1/2},t)$
the parameter $\alpha_\ell$ is defined so that $\tau^{\alpha_\ell}=2^\ell\tau^{1/2}\le 2^L\tau^{1/2}\le 2^{-K}\tau^2<\tau^{2/3}$, therefore $\alpha_\ell<\frac23$ as needed in Proposition  \ref{L:Sbound} for $d=3$.  We also have some restriction in $\tau$, they are (notice that  $\alpha_\ell=\frac12+\ell\frac{\log2}{\log\tau}$)
\[\frac{3}{4}\Bigl(\frac12+\ell\frac{\log2}{\log\tau}\Bigr)\log\tau\ge\log 3,\quad 
\Bigl\{\frac14+\frac16-\frac12\Bigl(\frac12+\ell\frac{\log2}{\log\tau}\Bigr)\Bigr\}\log\tau\ge\log(\widehat{C_3}/{M}),\]
\[\Bigl\{\frac16+\Bigl(\frac{1}{4}-\frac12\Bigr)\Bigl(\frac12+\ell\frac{\log2}{\log\tau}\Bigr)\Bigr\}\log\tau\ge \log(\widehat{A_3}/{M}).\]
for $0\le \ell\le L$. But we only have to check the first one for $\ell=0$ and the others for $\ell=L$. Equivalently, we need $\tau$ to satisfy
\[\tau\ge 3^{8/3},\quad \tau\ge 2^{3L}(\widehat{C_3}/{M})^6,\quad \tau\ge
2^{6L}(\widehat{A_3}/{M})^{24}.\]
By the definitions of $L$ and $K$, we have 
\[2^L<2^{-K}\tau^{3/2}<\tau^{2/3-2}\tau^{3/2}=\tau^{1/6}.\]
Therefore, we will take $M=\widehat{A_3}>\widehat{B_3}$, so that the last condition is automatically satisfied, and we need
\[\tau\ge 3^{8/3}=18.7208\dots, \quad \tau\ge(\widehat{C_3}/\widehat{A_3})^{12}=19.2088\]
(taking the exact values of the constant given in equations (21) and (15) of \cite{A181}).
Assuming this, we will have
\[S_1\le \sum_{\ell=0}^L(2^\ell\tau^{1/2})^{-\sigma}S(2^\ell\tau^{1/2},t)\le 
\widehat{A_3}\sum_{\ell=0}^L(2^\ell\tau^{1/2})^{-\sigma}\tau^{\frac16+\frac12\alpha_\ell}=
\widehat{A_3}\tau^{\frac16}\sum_{\ell=0}^L(2^\ell\tau^{1/2})^{\frac12-\sigma}.\]
That is, 
\begin{equation}\label{E:partial2}
S_1\le \widehat{A_3}\tau^{\frac16+\frac12(\frac12-\sigma)}\sum_{\ell=0}^L2^{\ell(\frac12-\sigma)}=\frac{\widehat{A_3}2^{\sigma-\frac12}}{2^{\sigma-\frac12}-1}\tau^{\frac12(\frac56-\sigma)}.
\qedhere\end{equation}
\end{proof}

\begin{proposition}\label{P:vanderCorput}
For $\tau=t/2\pi> 20$ and $\sigma\ge1$ we have 
\begin{equation}\label{E:mainbound}
\begin{aligned}
\Bigl|&\zeta(s)-\sum_{n\le \tau^{1/2}}\frac{1}{n^s}\Bigr|\le \\
& 
\frac{\widehat{A_3}2^{\sigma-\frac12}}{2^{\sigma-\frac12}-1}\tau^{\frac12(\frac56-\sigma)}+
\frac{2^{\sigma+1}A}{1-2^{-\sigma}}\tau^{\frac23(\frac34-\sigma)}+
\frac{2^{\sigma+3}A}{3\log 2}\tau^{\frac23(\frac14-\sigma)}\log\tau+
2^{\sigma+1}A\tau^{\frac23(\frac14-\sigma)}+8\tau^{2(\frac12-\sigma)}.
\end{aligned}
\end{equation}
\end{proposition}

\begin{proof}
We proceed as in the proof of Proposition \ref{P:van1} to get \eqref{E:partial1} and \eqref{E:partial2}. In our case with $\sigma\ge1$ and $2^K\le 2\tau^{4/3}$ we have 
\begin{align*}
S_2&\le 2A\tau^{\frac12-2\sigma}\sum_{k=1}^K2^{k\sigma}+4A\tau^{\frac32-2\sigma}\sum_{k=1}^K 2^{k(\sigma-1)}\le 2A\tau^{\frac12-2\sigma}\frac{2^{K\sigma}-1}{1-2^{-\sigma}}+
4A\tau^{\frac32-2\sigma} K 2^{K(\sigma-1)}\\
&\le 2A\tau^{\frac12-2\sigma} \frac{2^\sigma}{1-2^{-\sigma}}\tau^{\frac43\sigma}+
4A\tau^{\frac32-2\sigma}\Bigl(1+\frac{4\log\tau}{3\log 2}\Bigr)2^{\sigma-1}\tau^{\frac43(\sigma-1)}\\
&=\frac{2^{\sigma+1}A}{1-2^{-\sigma}}\tau^{\frac23(\frac34-\sigma)}+2^{\sigma+1}A
\Bigl(1+\frac{4\log\tau}{3\log 2}\Bigr)\tau^{\frac23(\frac14-\sigma)}\\
&=\frac{2^{\sigma+1}A}{1-2^{-\sigma}}\tau^{\frac23(\frac34-\sigma)}+
\frac{2^{\sigma+3}A}{3\log 2}\tau^{\frac23(\frac14-\sigma)}\log\tau+
2^{\sigma+1}A\tau^{\frac23(\frac14-\sigma)}.
\end{align*}
The sum $S_1$ is bounded exactly as in Proposition \ref{P:van1} assuming that $\tau\ge20$, and this finishes the proof.
\end{proof}
\begin{corollary}\label{C:end}
For $1\le\sigma\le 3/2$ and $\tau>20$ we have 
\begin{equation}
\Bigl|\zeta(s)-\sum_{n\le \tau^{1/2}}\frac{1}{n^s}\Bigr|\le 
39.209\;\tau^{-\frac{1}{12}}+24.447 \tau^{-\frac16}
+30.400\tau^{-\frac12}\log\tau+15.804 \tau^{-1/2}+8\tau^{-1}.\end{equation}
\end{corollary}
\begin{proof} 
We have substituted in \eqref{E:mainbound} the coefficients depending on $\sigma$ by its maximum on the interval $[1,3/2]$, obtaining 
\begin{multline*}
\Bigl|\zeta(s)-\sum_{n\le \tau^{1/2}}\frac{1}{n^s}\Bigr|\\ \le 
\frac{2\widehat{A_3}}{2-\sqrt{2}}\;\tau^{\frac12(\frac56-\sigma)}+
\frac{16A}{2\sqrt{2}-1}\tau^{\frac23(\frac34-\sigma)}
+\frac{16A\sqrt{2}}{3\log2}\tau^{\frac23(\frac14-\sigma)}\log\tau
+4\sqrt{2}A\tau^{\frac23(\frac14-\sigma)}
+8\tau^{2(\frac12-\sigma)}
\end{multline*}
Then we notice that all the functions on the right side are decreasing functions of $\sigma$ for $\sigma>1$, so the inequalities are true also putting on the right $\sigma=1$
Computing the numerical values of these coefficients and substituting the upper bounds we obtain 
\[\Bigl|\zeta(s)-\sum_{n\le \tau^{1/2}}\frac{1}{n^s}\Bigr|\le 
39.209\;\tau^{-\frac{1}{12}}+24.447 \tau^{-\frac16}
+30.400\tau^{-\frac12}\log\tau+15.804 \tau^{-1/2}+8\tau^{-1}.\qedhere
\]
\end{proof}

\section{End of proof} 

Finally, we need a lower value of $\zeta(\sigma+it)$. We use the one in Carneiro et al. \cite{HV} which is based on Trudgian \cite{Tr2} and improves on the previous values of  T. Trudgian \cite{Tr2} and  M.~Mossinghoff and T.~Trudgian \cite{MT}.

\begin{theorem}\label{T:TM}
For $\sigma\ge1$ and $t\ge 500$, we have
\begin{equation}
|\zeta(\sigma+it)|^{-1}\le 42.9\log t, \qquad \sigma\ge 1-\frac{1}{12\log t},\quad t\ge 132.16
\end{equation}
\end{theorem}

Carneiro et al. only prove the inequality for $\sigma=1$, but the same proof with slight modifications gives the above result.

\begin{proof}[Proof of Theorem \ref{T:main}]
We prove Theorem \ref{T:main} applying Rouche's Theorem to the rectangle $R=[1,3/2]\times[t_0,T]$
and the functions $\Rzeta(s)$ and $\zeta(s)$. Since we know that $\zeta(s)$ have no zeros in  $R$ it follows that $\Rzeta(s)$ do not vanish there. To this end, we must prove the inequality 
$|\Rzeta(s)-\zeta(s)|<|\zeta(s)|$ for any $s$ in the  boundary of $R$.

We proceed in the following way. For any $s$ in the boundary of $R$ with $\sigma\le 3/2$ and $t\ge 2707$ we have
\[|\Rzeta(s)-\zeta(s)|\le \Bigl|\Rzeta(s)-\sum_{n\le \tau^{1/2}}\frac{1}{n^s}\Bigr|+
\Bigl|\sum_{n\le \tau^{1/2}}\frac{1}{n^s}-\zeta(s)\Bigr|.\]
For the first term, we still apply Proposition \ref{P:RzetaApp},  to the second term we apply Corollary \ref{C:end}. All this machinery gives us 
\[|\Rzeta(s)-\zeta(s)|\le \tau^{-\sigma/2}+39.209\;\tau^{-\frac{1}{12}}+24.447 \tau^{-\frac16}
+30.400\tau^{-\frac12}\log\tau+15.804 \tau^{-1/2}+8\tau^{-1}.\]
Carneiro et al. inequality in Theorem \ref{T:TM} gives us for $s\in\partial R$
\[\frac{1}{42.9\log t}\le |\zeta(s)|.\]
It is clear that for $t$ big enough we may complete the reasoning showing that
\[\tau^{-\sigma/2}+39.209\;\tau^{-\frac{1}{12}}+24.447 \tau^{-\frac16}
+30.400\tau^{-\frac12}\log\tau+15.804 \tau^{-1/2}+8\tau^{-1}<\frac{1}{42.9\log t}.\]
If the inequality is true for a given $t$ and $\sigma=1$, then it is true for all $1\le\sigma\le3/2$. 
In fact, this inequality is true for 
\[\tau\ge 6.24072032490448651663628063807879324939223120097\times 10^{64}.\]
This corresponds to $t>t_0=2\pi\tau_0=3.9211\dots10^{65}$.
\end{proof}

\end{document}